\documentclass{conm-p-l}
\newtheorem{theorem}{Theorem}[section]

\newtheorem{Prop}[theorem]{Proposition}
\newtheorem{Cor}[theorem]{Corollary}

\theoremstyle{definition}

\newtheorem{example}[theorem]{Example}

\newtheorem{conj}{Conjecture}
\newtheorem{problem}{Problem}

\theoremstyle{remark}

\numberwithin{equation}{section}

\newcommand{\M}{{\mathcal{M}}}
\newcommand{\Z}{\mathbb{Z}}
\newcommand{\g}{{\mathfrak{g}}}
\newcommand{\good}{{\mathrm{good}}}
\renewcommand{\c}{{\mathrm{c}}}

\newcommand{\vphi}{\varphi}
\newcommand{\eps}{\varepsilon}
\newcommand{\wt}{\operatorname{wt}}
\newcommand{\te}{{\widetilde{e}}}
\newcommand{\tf}{{\widetilde{f}}}
\newcommand{\ba}{\begin{array}}
\newcommand{\ea}{\end{array}}
\newcommand{\eq}[1]{$$\ba{#1}}
\newcommand{\eneq}{\ea$$}
\newcommand{\la}{\lambda}
\newcommand{\La}{\Lambda}
\newcommand{\gl}{{\mathfrak{gl}}}
\newcommand{\eqt}{\begin{eqnarray}&&}
\newcommand{\eneqt}{\end{eqnarray}}
\newcommand{\ran}{\rangle}
\newcommand{\lan}{\langle}
\newcommand{\hs}{\hspace*}
\newcommand{\cl}{\colon}
\newcommand{\isoto}{\xrightarrow{\,\,\sim\,\,}}
\renewcommand{\P}{{\mathcal{P}}}
\newcommand{\into}{\hookrightarrow}
\newcommand{\N}{\mathbb{N}}
\newcommand{\gel}{{\mathfrak{l}}}
\newcommand{\Hom}{\operatorname{Hom}}

\newenvironment{tenumerate}{
  \begin{enumerate}
  
  }{\end{enumerate}}
\newenvironment{anumerate}{
  \begin{enumerate}
  
  }{\end{enumerate}}

\newcommand{\bnum}{\begin{tenumerate}}
\newcommand{\enum}{\end{tenumerate}}

\newcommand{\banum}{\begin{anumerate}}
\newcommand{\eanum}{\end{anumerate}}

\begin{document}

\title{Realizations of Crystals}
\author{Masaki Kashiwara}
\address{Research Institute for Mathematical Sciences,
Kyoto University, Kyoto, Japan}
\dedicatory{Dedicated to Professor Ryoshi Hotta on his sixtieth birthday}


\maketitle

\section{Introduction}
Once we know the crystal associated with representations,
we can read for example the branching rule, decomposition into irreducible
components of tensor products, etc.
Hence it is important to give the explicit crystal structure of
representations.
Several constructions are already known:
\begin{enumerate}
\item the description by using the elementary crystals
$B_i$ (\cite{K6}),
\item Littelmann's path realizations (\cite{Lit,Lit2}),
\item Young tableaux in $\gl_n$-case, and their variants 
in the classical Lie algebra case (\cite{KN}).
\item the realization using perfect crystals in the affine case (\cite{(KMN)^2}),
\end{enumerate}

In this paper, we shall add two more realizations of crystals.

\begin{enumerate}
\item[(5)] the monome realization,
\item[(6)] the lattice realization.
\end{enumerate}
The realization (5) is the one due to Hiraku Nakajima,
and its variation.

\section{Preliminary}
Let $P$ be a weight lattice,
$\{\alpha_i\in P;i\in I\}$ the set of simple roots, and
$\{h_i\in P^*;i\in I\}$ the set of simple co-roots.
Let $P_+$ be the set of dominant weights,
$P_+:=\{\la\in P;\mbox{$\lan h_i,\la\ran\ge0$ for any $i$}\}$.
Let $\Lambda_i\in P$ be the fundamental weight:
$\lan h_j,\Lambda_i\ran=\delta_{ij}$.

We shall not recall here the notion of crystals
(see e.g. \cite{K10}).
We shall call $I$-crystal when we want to emphasize the index set $I$
of simple roots.

For a crystal $B$ let us denote by $B^\vee$ the crystal obtained 
from $B$ by reversing the arrows :
$B^\vee{:=}\{b^\vee;b\in B\}$ with
$\wt(b^\vee)=-\wt(b)$, $\eps_i(b^\vee)=\vphi_i(b)$,
$\vphi_i(b^\vee)=\eps_i(b)$, $\te_i(b^\vee)=(\tf_ib)^\vee$,
$\tf_i(b^\vee)=(\te_ib)^\vee$.

For a dominant integral weight $\la\in P_+$, let us denote by $B(\la)$
the crystal associated with the irreducible representation $V(\la)$
with highest weight $\la$.
Then $B(-\la):=B(\la)^\vee$ is the crystal
associated with the irreducible representation $V(-\la)$
with lowest weight $-\la$.
Let us denote by $B(\infty)$ the crystal associated with the negative part
$U_q^-(\g)$ of the quantized algebra $U_q(\g)$.
Then $B(-\infty){:=}B(\infty)^\vee$ is the crystal
associated with the positive part $U_q^+(\g)$.
For $\la\in P$, we denote by $T_\la$ the crystal
$\{t_\la\}$ with $\wt(t_\la)=\la$,
$\eps_i(t_\la)=\vphi_i(t_\la)=-\infty$ and
$\te_i(t_\la)=\tf_i(t_\la)=0$.
Then there is a strict embedding
$B(\la)\into B(\infty)\otimes T_\la$.

For an $I$-crystal $B$ and a subset $J$ of $I$, let us denote by $\Psi_J(B)$
the $J$-crystal $B$ with the induced crystal structure
from $B$.

A crystal $B$ is called {\em semi-normal} if 
for each $i\in I$ the $\{i\}$-crystal $\Psi_{\{i\}}(B)$
is a crystal associated with an integrable module.
This is equivalent to saying that
$\eps_i(b)=\max\{n\in \N\,;\,\te_i^nb\not=0\}$
and $\vphi_i(b)=\max\{n\in \N\,;\,\tf_i^nb\not=0\}$ for any $b\in B$ and
$i\in I$.
A crystal $B$ is called {\em normal} if for any subset $J$ of $I$
such that $\{\alpha_i;i\in J\}$ is a root system for finite-dimensional 
semisimple Lie algebra,
$\Psi_J(B)$ is the crystal associated with an integrable
$U_q(\g_J)$-module,
where $U_q(\g_J)$ is the associated quantum group with $\{\alpha_i;i\in J\}$.

The crystal $B_i$ is $\{b_i(n);n\in\Z\}$
with $\wt(b_i(n))=n\alpha_i$,
$$
\ba{ll}
\eps_j(b_i(n))=
\left\{\ba{ll}
-n&\mbox{for $j=i$}\\
-\infty&\mbox{otherwise}
\ea
\right.
&
\vphi_j(b_i(n))=
\left\{\ba{ll}
n&\mbox{for $j=i$}\\
-\infty&\mbox{otherwise}
\ea
\right.\\[5pt]
\te_j(b_i(n))=
\left\{\ba{ll}
b_i(n+1)&\mbox{for $j=i$}\\
0&\mbox{otherwise}
\ea
\right.&
\tf_j(b_i(n))=
\left\{\ba{ll}
b_i(n-1)&\mbox{for $j=i$}\\
0&\mbox{otherwise}
\ea
\right.
\ea
$$
Let $\{i_n\}_{n=1,2,\ldots}$ be a sequence in $I$ such that
every element in $I$ appears infinitely many times in this sequence.
Then $B(\infty)$ is strictly embedded into the crystal
$\cdots\cdots\otimes B_{i_2}\otimes B_{i_1}$.
This is the realization (1) in the introduction.

\medskip
The following theorem is frequently used in this paper.
\begin{theorem}[\cite{{(KMN)^2},SMF}]
Let $B$ be a crystal satisfying the condition:
for any subset $J$ of $I$ with at most two elements,
any connected component of $\Psi_J(B)$
containing a highest weight vector
{\rm(}i.e.\ a vector annihilated by any $\te_i (i\in J)${\rm)}
is a crystal isomorphic to the crystal associated with
an integrable highest weight $U_q(\g_J)$-module.

Then
any connected component of $B$ containing a highest weight vector
is isomorphic to $B(\la)$ for a dominant integral weight $\la$.
\end{theorem}

The following variation is proved similarly.
\begin{theorem}
Let $B$ be a crystal satisfying the condition:
for any subset $J$ of $I$ with at most two elements,
any connected component of $\Psi_J(B)$
containing a highest weight vector $b$
{\rm(}i.e.\ a vector annihilated by any $\te_i (i\in J)${\rm)}
is a crystal isomorphic to the $J$-crystal $B_J(\infty)\otimes T_{\wt(b)}$.
Here $B_J(\infty)$ is the crystal associated with
the negative part $U_q^-(\g_J)$ of $U_q(\g_J)$.

Then
any connected component of $B$ containing a highest weight vector $b$
is isomorphic to $B(\infty)\otimes T_{\wt(b)}$ for a dominant integral weight $\la$.
\end{theorem}

\section{Mon\^ome realization}

The following realization is due to Hiraku Nakajima (\cite{Hi}).

Let $\M$ be the set of mon\^omes in the variables
$Y_i(n)$ ($i\in I$, $n\in\Z$):
$$\M{:=}\{\prod\limits_{i\in I,n\in\Z}Y_i(n)^{y_i(n)}; \mbox{
$y_i(n)\in\Z$ vanish except finitely many $(i,n)$}\}.$$
We set 
$A_i(n)=Y_i(n-1)Y_i(n+1)\prod_{k\not=i}Y_k(n)^{\langle h_k,\alpha_i\rangle}$.
We assign to $Y_i(n)$ the weight $\Lambda_i$,
so that $A_i(n)$ has the weight $\alpha_i$.
We define the crystal structure on $\M$ as follows.
For a mon\^ome $M=\prod\limits_{i\in I,n\in\Z}Y_i(n)^{y_i(n)}$,
we set
\eq{ll}
\wt(M)&=\sum_i(\sum_ny_i(n))\Lambda_i,\\[5pt]
\vphi_i(M)&=\max\{\sum_{k\le n}y_i(k);n\in\Z\},\\[5pt]
\eps_i(M)&=\max\{-\sum_{k\ge n}y_i(k);n\in\Z\}.
\eneq
We define
\eq{ll}
\tf_i(M)&=
\left\{\ba{ll}
0&\mbox{if $\vphi_i(M)=0$,}\\
A_i(n_f+1)^{-1}M&\mbox{if $\vphi_i(M)>0$,}
\ea\right.\\[10pt]
\te_i(M)&=
\left\{\ba{ll}
0&\mbox{if $\eps_i(M)=0$,}\\
A_i(n_e-1)M&\mbox{if $\eps_i(M)>0$.}
\ea\right.
\eneq
Here 
\eq{ll}
n_f&=\min\{n;\vphi_i(M)=\sum_{k\le n}y_i(k)\}\\[5pt]
    &=\min\{n;\eps_i(M)=-\sum_{k>n}y_i(k)\},\\[10pt]
n_e&=\max\{n;\vphi_i(M)=\sum_{k<n}y_i(k)\}\\[5pt]
    &=\max\{n;\eps_i(M)=-\sum_{k\ge n}y_i(k)\}.
\eneq
Note that $y_i(n_f)>0$, $y_i(n_f+1)\le0$ and $y_i(n_e)<0$, $y_i(n_e-1)\ge0$.

\medskip
\emph{It is in fact not true that $\M$ is a crystal}
(see Example \ref{ex:bad}).
All the axioms for crystals hold but
$b'=\tf_ib\Longleftrightarrow \te_ib'=b$.

\begin{conj}
For a product $M$ of positive powers of $Y_i(n)$,
the connected component
containing $M$
is isomorphic to $B(\wt(M))$.
\end{conj}

This is checked in numerous examples by computer.
The following fact is observed
by the computer experiment
in the connected component $B$ of the product of positive powers
 :
for any $i$ and $n$, $y_i(n)>0$, $y_i(n+1)<0$ can
never happen at once.
This property implies
$\tf_i\te_iM=M$ and $\te_i\tf_iM=M$ if they are not zero.

\medskip

In fact we have
\begin{Prop}
Let $\M_\good$ be a subset of $\M$ satisfying the following properties.
\begin{itemize}
\item[(i)]
For $M\in\M_\good$ and $i$ $n$,  $y_i(n)>0$ implies $y_i(n+1)\ge0$.
\item[(ii)]
$\M_\good$ is stable by the $\te_i$'s and the $\tf_i$'s.
\end{itemize}
Then $\M_\good$ has good properties.
Namely it is a normal crystal,
and the crystal generated by a highest weight vector of weight $\lambda$
is isomorphic to $B(\lambda)$.
\end{Prop}
Since the proof is similar to the proof of
Theorem \ref{th:M_c}, we omit this.

\bigskip
\begin{problem}
Find a subcrystal $\M_\good$ of $\M$
which has the properties above
(or modify the definition above so that $\M$ is really a crystal).
\end{problem}

\bigskip
\begin{example} For $\g=\mathfrak{sl}_3$,
$$Y_1(0)\xrightarrow{1}Y_1(2)^{-1}Y_2(1)\xrightarrow{2}Y_2(3)^{-1}$$
is the crystal isomorphic to $B(\Lambda_1)$.
\end{example}

\begin{example}\label{ex:bad}
For $\M$, all the axioms for crystals hold except
$$b'=\tf_ib\Longleftrightarrow \te_ib'=b.$$
For example, for $\g=\mathfrak{sl}_2$.
$$Y(1)Y(2)^{-1}\xrightarrow{\tf_1}Y(2)^{-1}Y(3)^{-1}$$
and the chain in $\M_\good$
$$Y(0)Y(1)\xrightarrow{\tf_1}Y(0)Y(3)^{-1}
\xrightarrow{\tf_1}Y(2)^{-1}Y(3)^{-1}.$$
So $Y(1)Y(2)^{-1}\not\in\M_\good$, or $\tf_1(Y(1)Y(2)^{-1})=0$
(with another modified rule).
\end{example}

%

\section{Variations}

We shall define another structure of crystal on $\M$.
Let $\c=\bigl(\c_{ij}\bigr)_{i\neq j\in I}$
be a set of integers such that
$$\c_{ij}+\c_{ji}=1.$$

We set
$$A_i(n)=Y_i(n)Y_i(n+1)
\prod\limits_{j\neq i}Y_j(n+\c_{ji})^{\lan h_j,\alpha_i\ran}.
$$

For a mon\^ome $M=\prod\limits_{i\in I,n\in\Z}Y_i(n)^{y_i(n)}$,
we set
\eq{ll}
\wt(M)&=\sum_i(\sum_ny_i(n))\Lambda_i,\\[5pt]
\vphi_i(M)&=\max\{\sum_{k\le n}y_i(k);n\in\Z\},\\[5pt]
\eps_i(M)&=\max\{-\sum_{k>n}y_i(k);n\in\Z\}.
\eneq
We define
\eq{ll}
\tf_i(M)&=
\left\{\ba{ll}
0&\mbox{if $\vphi_i(M)=0$,}\\
A_i(n_f)^{-1}M&\mbox{if $\vphi_i(M)>0$,}
\ea\right.\\[10pt]
\te_i(M)&=
\left\{\ba{ll}
0&\mbox{if $\eps_i(M)=0$,}\\
A_i(n_e)M&\mbox{if $\eps_i(M)>0$.}
\ea\right.
\eneq
Here 
\eq{ll}
n_f&=\min\{n;\vphi_i(M)=\sum_{k\le n}y_i(k)\}\\[5pt]
    &=\min\{n;\eps_i(M)=-\sum_{k>n}y_i(k)\},\\[10pt]
n_e&=\max\{n;\vphi_i(M)=\sum_{k\le n}y_i(k)\}\\[5pt]
    &=\max\{n;\eps_i(M)=-\sum_{k>n}y_i(k)\}.
\eneq
Note that $y_i(n_f)>0$, $y_i(n_f+1)\le0$ and $y_i(n_e+1)<0$, $y_i(n_e)\ge0$.

\noindent
Let us denote by $\M_\c$ the crystal $\M$ thus defined.

In this case, one can prove easily the following peroposition.
\begin{Prop}
$\M_\c$ is a semi-normal crystal.
\end{Prop}

For a family of integers $(m_i)_{i\in I}$,
let us set $\c'=(c'_{ij})_{i\neq j\in I}$,
by $c'_{ij}=c_{ij}+m_i-m_j$.
Then the map $Y_i(n)\mapsto Y_i(n+m_i)$
gives a crystal isomorphism
$\M_\c\isoto \M_{\c'}$.
Hence if the Dynkin diagram has no loop,
the isomorphism class of the crystal
$\M_\c$ does not depend on the choice of $c$.

Let $\psi\cl \M\to\M$ be the map defined by
$Y_i(n)\mapsto Y_i(-n)^{-1}$.
Then one has

\begin{Prop}
The map $\psi$ give an isomorphism of crystals:
${\M_\c}^\vee\isoto\M_{{\mathrm{c'}}}$
where ${\mathrm{c'}}_{ij}=\c_{ji}$.
Here ${\M_\c}^\vee$ is the crystal obtained from $\M_\c$
by reversing the direction of arrows.
\end{Prop}

\begin{theorem}\label{th:M_c}
For a highest wieght vector $M\in\M_c$,
the connected component of $\M_\c$
containing $M$
is isomorphic to $B(\wt(M))$.
\end{theorem}
\begin{proof}
It is enough to show for a subset
$J$ of $I$ with cardinality $2$,
any connected component
$B$ of $\M_\c$ containing
a highest weight vector
is isomorphic to the crystal
associated with an integrable highest wieight
$U_q(\g_J)$-module. 
Assume $J=\{1,2\}$.
Then we may assume that $c_{1,2}=0$ and $c_{2,1}=1$.
Then one has
$$A_1(n)=Y_1(n)Y_1(n+1)Y_2(n+1)^{\lan h_2,\alpha_1\ran}\times\,\cdots$$
and
$$A_2(n)=Y_2(n)Y_2(n+1)Y_1(n)^{\lan h_1,\alpha_2\ran}\times\,\cdots.$$

Take $\la(n)=\sum_{i}\la_i(n)\La_i\in P$ ($n\in \Z$),
which vanish except for finitely many $n$.
Set $K_n=B_1\otimes B_2\otimes T_{\la(n)}$.

$$K:=\cdots\otimes K_{1}\otimes K_0\otimes K_{-1}\otimes \cdots.$$
We define the map $\Phi$ from $K$ to $\M_\c$ by
$$\ba{l}
K\ni b{:=}\mathop\bigotimes\limits_n
\bigl(b_1(z_1(n))\otimes b_2(z_2(n))\otimes t_{\la(n)}\bigr)\\
\hs{60pt}\mathop{\longmapsto}\limits^\Phi
M{:=}\prod_{i,n}Y_i(n)^{\la_i(n)}A_1(n)^{z_1(n)}A_2(n)^{z_2(n)}\in\M_\c.\ea$$
Hence
$$\ba{ll}
y_i(n)&=\sum_n\Bigl(\la_i(n)+z_i(n)+z_i(n-1)+
\sum_{\nu\neq i}\lan h_i,\alpha_\nu\ran z_\nu(n-c_{i,\nu})\Bigr)\\[5pt]
&=\left\{
\ba{ll}
\la_1(n)+z_1(n)+z_1(n-1)+
\lan h_1,\alpha_2\ran z_2(n)&\quad\mbox{for $i=1$,}\\[3pt]
\la_2(n)+z_2(n)+z_2(n-1)+
\lan h_2,\alpha_1\ran z_1(n-1)&\quad\mbox{for $i=2$.}
\ea\right.\ea$$
One has
$$\eps_1(b)
=\max\{-z_1(n)
-\sum_{k>n}\bigl(2z_1(k)+\lan h_1,\alpha_2\ran z_2(k)+\la_1(k)\bigr)
\,;\,n\in \Z\}.$$
On the other hand
$$
\ba{ll}\eps_1(M)&=\max\{-\sum_{k>n}y_1(k)\,;n\in \Z\}\\[3pt]
&=\max\bigl\{-\sum_{k>n}\bigl(\la_1(k)+z_1(k)+z_1(k-1)+
\lan h_1,\alpha_2\ran z_2(k)\bigr)
;n\in\Z\bigr\}.\ea$$
Hence $\eps_1(b)$ and $\eps_1(M)$ are equal.
One has
$$
\ba{ll}\eps_2(b)
&=\max\{-z_2(n)-\lan h_2,\alpha_1\ran z_1(n)\\[3pt]
&\hs{60pt}-\sum_{k>n}\bigl(2z_2(k)+\lan h_2,\alpha_1\ran z_1(k)+\la_2(k)\Bigr)
\,;\,n\in \Z\},\ea$$
and
$$
\ba{ll}
\eps_2(M)&=\max\{-\sum_{k>n}y_2(k)\,;n\in \Z\}\\[3pt]
&=\max\{-\sum_{k>n}\bigl(\la_2(k)+z_2(k)+z_2(k-1)+
\lan h_2,\alpha_1\ran z_1(k-1)\bigr)
;n\in\Z\},\ea$$
and hence $\eps_2(b)$ and $\eps_2(M)$ are equal.
It is easy to see the action of $\tf_i$ and $\te_i$ commute with $\Phi$.
Hence $\Phi$ is a crystal morphism.

On the other hand, one has.
$T_\la\otimes B_i\simeq B_i\otimes T_{s_i\la}$ and 
$T_\la\otimes T_\mu\simeq T_{\la+\mu}$.
Hence, one has an isomorphism
$$K\simeq \P\otimes T_\la\otimes \P^-$$
for some $\la$.
Here $B=B_1\otimes B_2$, $\P=\cdots \otimes B\otimes B$ and  
$\P^-=B\otimes B\otimes \cdots$.
One knows
$\P=\bigsqcup\limits_{b\in \P^h}B(\infty)\otimes T_{\wt(b)}$
and $\P^-=\bigsqcup\limits_{b\in (\P^-)^\ell}B(-\infty)\otimes T_{\wt(b)}$,
where
$\P^h=\{b\in \P;\mbox{$\te_ib=0$ for any $i$}\}$
and $(\P^-)^\ell=\{b\in (\P^-)^\ell;\mbox{$\tf_ib=0$ for any $i$}\}$.
Hence $K$ is a disjoint union of $B(\infty)\otimes T_\la\otimes B(-\infty)
\subset B(\tilde U_q(\g))$.
Hence the connected component containing a highest weight vector is isomorphic to $B(\la)$ for some $\la\in P_+$.
\end{proof}

\begin{Cor}
$\M_\c$ is a normal crystal.
\end{Cor}
The condition
$\c_{ij}+\c_{ji}=1$
might be relaxed in the non simply-laced case.

For example, it seems that
Theorem \ref{th:M_c} holds when $\c_{ij}+\c_{ji}\ge1$
in the $A_1^{(1)}$-case,
checked by computer experiment.

\begin{problem}
What is the condition for $\c$ in order that
the above theorem holds?
\end{problem}

\begin{problem}
What is the connected component
of $\M_\c$ containing a given highest weight vector?
\end{problem}

\section{Lattice realization}

Let $Q$ be the root lattice $\sum_{i\in I}\Z\alpha_i$.
Let us take a family $\gel=\{\ell_i\}_{i\in I}$
of elements in $Q^*{:=}\Hom(Q,\Z)$, 
which satisfies the following property:
\begin{eqnarray}
\text{$\ell_i(\alpha_i)=-1$ for every $i\in I$.}
\end{eqnarray}

Let us define the crystal $B_\gel=\{v(x);x\in Q\}$ as follows.
One has
$\wt(v(x))=x$, $\eps_i(v(x))=\ell_i(x)$, 
$\vphi_i(b)=\lan h_i,\wt(b)\ran+\eps_i(b)$,
$\te_i(v(x))=v(x+\alpha_i)$ and $\tf_i(v(x))=v(x-\alpha_i)$.

\begin{theorem}
Assume further that $\gel$ satisfies the following condition:
\eqt\label{eq:ell}
\ba{l}
\mbox{for every pair $(i,j)$ with $i\neq j$,
one has either}\hspace{100pt}\\[5pt]
\qquad\mbox{{\rm (i)}
$\ell_i(\alpha_j)=-\lan h_i,\alpha_j\ran$ and $\ell_j(\alpha_i)=0$, or}
\\[3pt]
\qquad\mbox{{\rm (ii)}
$\ell_i(\alpha_j)=0$ and $\ell_j(\alpha_i)=-\lan h_j,\alpha_i\ran$.}
\ea
\eneqt
Then there is a strict embedding
$B(\infty)\into B(\infty)\otimes B_\gel$ sending
$u_\infty$ to $u_\infty\otimes v(0)$.
Hence $B(\infty)$ is embedded into
$\cdots\otimes B_\gel\otimes B_\gel$.
\end{theorem}

\begin{proof}
We may assume $I=\{1,2\}$.
Assume that $\ell_2(\alpha_1)=-\lan h_2,\alpha_1\ran$
and $\ell_1(\alpha_2)=0$.
Then $B_\gel\simeq B_1\otimes B_2$, and
the result follows from the realization (1).
\end{proof}

\begin{example}
For $\g=A^{(1)}_n$ ($n\ge2$) with
$I=\Z/(n+1)\Z$, if we take $\ell_i(\alpha_j)=-\delta_{i,j}+\delta_{i+1,j}$,
then $B_\gel$ is the crystal associated with
a coherent family of perfect crystals (\cite{KKM}).
Hence the above result is proved in loco.cito.
\end{example}
\begin{example}
For $\g=A^{(1)}_1$ with $I=\{0,1\}$, if we take 
$$\ell_i(\alpha_j)=\left\{\ba{rl}-1&\mbox{if $i=j$,}\\
1&\mbox{for $i\neq j$,}\ea\right.$$
then $B_\gel$ is the crystal associated with
a coherent family of perfect crystals (\cite{KKM}).
Hence $B(\infty)$ is embedded into
$\cdots\otimes B_\gel\otimes B_\gel$,
although $\gel$ does not satisfy the condition \eqref{eq:ell} above. 
\end{example}
\begin{problem}
What is the condition for $\gel$ in order that the
above theorem holds?
\end{problem}


\begin{thebibliography}{99}
\providecommand{\bysame}{\makebox[3em]{\hrulefill}\thinspace}
%
%
%
%
%
%
%
%
%
%
%
%
%
%
\bibitem{K6} M. Kashiwara,
\emph{The crystal bases and Littelmann's refined Demazure character formula},
Duke Math. J. $\mathbf{71}$ (1993), $\textrm{n}^\circ 3$, 839--858.
%
%

\bibitem{K10} \bysame,
\emph{On crystal bases},
Representations of Groups,
Proceedings of the 1994 Annual Seminar of the Canadian Math. Soc.
Banff Center, Banff, Alberta,
June 15--24,
(B.N. Allison and G.H. Cliff, eds), 
CMS Conference proceedings, {\bf 16} (1995) 155--197, Amer. Math. Soc., 
Providence, RI.

\bibitem{SMF}
\bysame,
\emph{Bases cristallines des groupes quantiques},
noted by Charles Cochet,
submitted to ``cours specialis\'e, SMF.



\bibitem{KKM}
S.-J. Kang, M.  Kashiwara and K. Misra,
\emph{Crystal bases of Verma modules for quantum affine Lie algebras},
Compositio Math. {\bf 92} (1994), no. 3,
299--325.

\bibitem{(KMN)^2}
M. Kashiwara, S.-J.Kang, K.Misra, T.Miwa, T.Nakashima and A. Nakayashiki,
Affine crystals and vertex models,
International Journal of Modern Physics A {\bf 7},
Suppl.1A(1992) 449--484,
Proceeding of the RIMS Research Project 1991 ``Infinite Analysis''. 

\bibitem{KN} M. Kashiwara and T. Nakashima,
\emph{Crystal graphs for representations of the $q$-analogue 
of classical Lie algebras},
J. of Algebra, $\mathbf{165}$ (1994), 295--345.

\bibitem{Lit} P. Littelmann,
\emph{A Littlewood-Richardson rule for symmetrisable Kac-Moody Lie algebras},
Invent. Math. $\mathbf{116}$ (1994), 329--346.

\bibitem{Lit2} \bysame,
\emph{Paths and root operators in representation theory},
Ann. of Math. (2) $\mathbf{142}$ (1995), no. 3, 499--525.
%
%
%
%

\bibitem{Hi} H. Nakajima,
\emph{$t$--analogs of $q$--characters of quantum affine algebras of type
$A_n$, $D_n$}, in the same volume.


\bibitem{Na} T. Nakashima,
\emph{Crystal bases and a generalization of Littlewood-Richardson rule for the classical Lie algebras},
Commun. Math. Phys. $\mathbf{154}$ (1993), 215--243.
%
%
\end{thebibliography}
\end{document}